\documentclass[reqno, a4paper, 12pt]{amsart}
\usepackage{amsmath, amssymb, eucal, amscd, amstext, enumerate}
\usepackage{amsfonts,amsthm}
\usepackage{mathrsfs}

\usepackage[numbers,sort&compress]{natbib}
\usepackage{color,hyperref}

\theoremstyle{plain}
\newtheorem{thm}{Theorem}[section]
\newtheorem{lem}[thm]{Lemma}

\newtheorem{defn}[thm]{Definition}
\newtheorem{rem}[thm]{Remark}
\newtheorem{rem-ntn}[thm]{Remark and Notation}

\newenvironment{prf}{{\noindent \textbf{Proof:}\ }}{\hfill $\Box$\\ \smallskip}

\numberwithin{equation}{section}

\newcommand{\mednoind}{{\medskip\noindent}}

\newcommand{\norm}[1]{\left\|#1\right\|}
\newcommand{\form}[1]{{\rm\bf{#1}}}

\newcommand{\KH}{\mathfrak{H}}

\newcommand{\BR}{\mathbb{R}}
\newcommand{\BN}{\mathbb{N}}

\newcommand{\BC}{\mathbb{C}}

\usepackage{color}

\begin{document}

\title[A new approach to inverse Sturm-Liouville problems II. The singular case.]
{A new approach to inverse Sturm-Liouville problems based on point interaction II. The singular case.}

\author[Min Zhao]{Min Zhao}
\address{Department of Mathematics, Shandong University at
         Weihai \\Weihai 264209, P.R. China}
\email{zhaomin215@mail.sdu.edu.cn}\

\author[Jiangang Qi]{Jiangang Qi}
\address{Department of Mathematics, Shandong University at
         Weihai \\Weihai 264209, P.R. China}
\email{qjg816@163.com}\

\author[Xiao Chen]{Xiao Chen$^\dag$}
\address{Department of Mathematics, Shandong University at
         Weihai \\Weihai 264209, P.R. China}
\email{chenxiao@sdu.edu.cn}\


\begin{abstract} In this paper, further to the point interaction method for inverse Sturm-Liouville problems on finite intervals firstly proposed in our previous work, we will continue to generalize this method to the inverse eigenvalue problems for singular Sturm-Liouville problems on the half real axis.

\vspace{04pt}


\end{abstract}

\begin{thanks}{$\dag$ the corresponding author, \emph{Email}: \texttt{chenxiao@sdu.edu.cn}.}
\end{thanks}

\maketitle

\section{Introduction and problem statement}\label{sec:intro}

In our previous paper \cite{ZQC2024}, we develop a new and explicit scheme to inverse eigenvalue problems for regular Sturm-Liouville equations under Dirichlet boundary.
More precisely, for a given Sturm-Liouville problem (abbreviated as  S-L problem) on a finite interval under Dirichlet boundary with the unknown integrable potential, by letting this problem interact with $\delta$-function potentials at every point in the interval of definition of solutions to this S-L problem, we obtain a family of perturbation problems, originating from {\bf moving point interaction models} in quantum mechanics (cf. \cite{Pos2007, AGHH1988}).
Then, only depending on the first eigenvalues of these perturbation problems, we define and study the {\bf first eigenvalue function}, by which the desired potential can be expressed explicitly and uniquely.

In this paper, we are devoted to generalize the method of inverse S-L problems proposed in \cite{ZQC2024} from regular S-L problems to singular cases on semi-infinity internals.

It is well known that, a S-L problem $({\bf \tilde{E}}_q)$ on the half real axis with regular endpoint 0 under the separation boundary condition is defined as:
\begin{equation}\label{eqn:general-probl}
\tau_q y:=\ -y''(x)+q(x)y(x)=\lambda y(x),\ x\in[0,+\infty),\ y(0)\sin\alpha-y'(0)\cos\alpha=0=\mathbb{B}y(+\infty),
\end{equation}
where $0\leq \alpha<\pi$, $q\in L^1_{loc}([0,+\infty),\BR)$, the complex number $\lambda$ is the spectral parameter and $\mathbb{B}y(+\infty)$ is the boundary condition at $+\infty$. Here $L^1_{loc}([0,+\infty),\BR)$ denotes the space of real value locally integrable functions on $[0,+\infty)$. See \cite[Section 10.4]{Zettl2005}.

As we all know, for the singular differential expression $\tau_q y$ defined on the interval $[0,+\infty)$, H.\ Weyl \cite{Weyl1910} divided it into two classes:
the limit-circle case where all solutions of $\tau_q y=\lambda y$ are in $L^2[0,+\infty)$ and the limit-point case where essentially only one solution is in $L^2[0,+\infty)$  for non-real values of the parameter.
Here $L^2[0,+\infty)$ denotes the space of complex value square-integrable functions on $[0,+\infty)$ equipped with the canonical norm $\norm{\cdot}_2$.

\medskip

The maximal domain of $\tau_q y$ is defined by
\begin{equation}\label{def-max-dom}
\mathcal{D}_{max}:=\left\{ y\in L^2[0,+\infty)\, :\, y,\, y'\in AC_{loc}[0,+\infty),\ \tau_q y\in L^2[0,+\infty) \right\},
\end{equation}
where $AC_{loc}[0,+\infty)$ is the space of locally absolutely continuous complex value functions on the interval $[0,+\infty)$ (see \cite[(10.2.2)]{Zettl2005}).
For any $y,z\in \mathcal{D}_{max}$, the limits
\begin{equation}\label{equ-lim}
[y,z](0)=\lim\limits_{x\rightarrow 0}[y,z](x),\ [y,z](+\infty)=\lim\limits_{x\rightarrow +\infty}[y,z](x)
\end{equation}
exist and are finite, where $[y,z]=y\overline{z}'-\overline{z}y'$ (see \cite[Lemma 10.2.3]{Zettl2005}).

If the endpoint $+\infty$ is regular, the limits \begin{equation}\label{equ-lim-inft}
y(+\infty)=\lim\limits_{x\rightarrow+\infty}y(x)\text{ and } y'(+\infty)=\lim\limits_{x\rightarrow+\infty}y'(x)
\end{equation}
exist and are finite for any function $y\in \mathcal{D}_{max}$ (see \cite[Lemma 10.4.3]{Zettl2005}). Then the self-adjoint domain of $\tau_q y$ is given by
\begin{equation}\label{def-oper-reg}
\mathcal{D}(\tilde{S}_q):=\left\{y\in \mathcal{D}_{max}\ :\ y(0)\sin\alpha-y'(0)\cos\alpha=0=y(+\infty)\sin\beta-y'(+\infty)\cos\beta,\ 0\leq\alpha,\ \beta<\pi \right\},
\end{equation}
where $\tilde{S}_q$ denotes the corresponding self-adjoint operator, see \cite[Theorem 10.4.3(1)]{Zettl2005}.

If the endpoint $+\infty$ is singular, then the limits \eqref{equ-lim-inft} do not exist.
Based on the separation boundary condition in \eqref{eqn:general-probl}, we can obtain the self-adjoint domain of $\tau_q y$ (see \cite[Theorem 10.4.4 and (10.4.60)]{Zettl2005}).
If $\tau_q y$ in $({\bf \tilde{E}}_q)$ is in the limit-point case at $+\infty$, the self-adjoint domain of $\tau_q y$ is given by
\begin{equation}\label{def-oper-lim-point}
\mathcal{D}(\tilde{S}_q):=\left\{y\in \mathcal{D}_{max}\ :\ y(0)\sin\alpha-y'(0)\cos\alpha=0 \right\},
\end{equation}
while, if $\tau_q y$ in $({\bf \tilde{E}}_q)$ is in the limit-circle case at $+\infty$, the self-adjoint domain of $\tau_q y$ is given by
\begin{equation}\label{def-oper-lim-circle}
\mathcal{D}(\tilde{S}_q):=\left\{y\in \mathcal{D}_{max}\ :\ y(0)\sin\alpha-y'(0)\cos\alpha=0=[y,f](+\infty)\sin\beta-[y,g](+\infty)\cos\beta \right\},
\end{equation}
where $\alpha$ is just the one in \eqref{eqn:general-probl}, $0<\beta\leq \pi$ and $f,g$ are (BC) basis at $+\infty$, that is, a pair of real-valued functions $f,\ g\in \mathcal{D}_{max}$ and satisfying $[f,g](+\infty)=1$ (see \cite[Definition 10.4.2]{Zettl2005}).

\bigskip

Similar to the regular case, for the unknown potential $q$, it is still impossible to recover $q$ uniquely by the first eigenvalue of the singular S-L problem.
Therefore, we will generalize the method in \cite{ZQC2024} to the singular case.

The key point of proving the main theorem is the existence of the {\bf non-oscillatory eigenvalue} of $({\bf \tilde{E}}_q)$, that is, the corresponding eigenfunction has finite number of zeros.
In fact, we just need to guarantee the existence of the {\bf principle eigenvalue} of $({\bf \tilde{E}}_q)$, that is, the corresponding eigenfunction has no zero on $(0, +\infty)$, i.e., this eigenfunction does not change sign on $(0, +\infty)$. Since our method in \cite{ZQC2024} mainly deals with the principle eigenvalue for regular S-L problems, we here only need to develop a method which can deal with the principle eigenvalue for singular S-L problems.
In other words, if there exists a non-oscillatory eigenvalue of $({\bf \tilde{E}}_q)$ such that the corresponding eigenfunction has $n$ zeros $x_1$, ..., $x_n$ which are ordered as $x_0=0<x_1<\cdots <x_n<+\infty$ on $(0,+\infty)$, then we can recover the unknown potential on interval $[x_i,x_{i+1}],\ i=0,\cdots,n-1$ using the method in \cite{ZQC2024}, and on $[x_n,+\infty)$ using the method to be proposed in the present paper.
For the benefit of the readers, the main steps of the latter method are now roughly described as follows.

In the present paper, we mainly consider the problem $({\bf \tilde{E}}_q)$ with $\alpha=\frac{\pi}{2}$, denoted by $({\bf E}_q)$, defined as
\begin{equation}\label{eqn:main-probl}
\tau_q y:=\ -y''(x)+q(x)y(x)=\lambda y(x),\ x\in[0,+\infty),\ y(0)=0=\mathbb{B}y(+\infty),
\end{equation}
where $q\in L^1_{loc}([0,+\infty),\BR)$ and $\mathbb{B}y(+\infty)$ are the same as in \eqref{eqn:general-probl}. Then we can also obtain the corresponding self-adjoint operator $S_q$ as above.

For the problem $({\bf  E}_q)$, we propose two hypotheses:
\begin{itemize}
 \item[] ${\bf (H_1)}:$ \emph{the potential $q$ is non-negative (equivalently $q$ is bounded below). Note that, in this case, $\tau_q y$ is in the limit-point case at $+\infty$ (see Lemma~\ref{lem:principle-eigenv-exist-case}(i)).}
 \item[] ${\bf (H_2)}:$ \emph{$\tau_q y$ is in the limit-circle case and non-oscillatory at $+\infty$.}
\end{itemize}

In the sequent section, we will see that each of the above hypotheses can ensure that $({\bf E}_q)$ has the first eigenvalue which is also  the principle eigenvalue (see Lemma~\ref{lem:principle-eigenv-exist-case}).

To begin with, inspired by $\delta$-point interaction in quantum mechanics, we embed $\delta$-function potential into the original system \eqref{eqn:main-probl}.

Denote by $\delta(x-t)$ the {\bf Dirac $\delta$-function} at $t\in(0,+\infty)$, which is defined by
\begin{equation}\label{eqn:delta-funct}
\delta(x-t)=
\begin{cases}
+\infty, &x=t,\\
          0, &x\neq t
\end{cases}
\ \text{and}\ \int_I\delta(x-t)\, {\rm d}x=1,\ \ \forall\, I\subset (0,+\infty)
\ \text{and}\ t\in I.
\end{equation}
And there holds $$\int_0^{+\infty} f(x)\delta(x-t)\,{\rm d}x=\int_{t-\epsilon}^{t+\epsilon} f(x)\delta(x-t)\,{\rm d}x=f(t),$$ for any continuous function $f$ and $0<\epsilon<t$.

The Dirac measure $\delta_t$ (also called the Dirac distribution, the point mass, or Heaviside function) at $t\in(0,+\infty)$, as a kind of discrete measure, is defined by
\begin{equation}\label{eqn:delta-mes}
\delta_t(x)=\left\{
\begin{array}{cl}
0, &  x\in[0,t), \\
1,  &  x\in[t,+\infty). \\
\end{array} \right.
\end{equation}
Note that, the Dirac $\delta$-function $\delta(x-t)$ defined by \eqref{eqn:delta-funct} is exactly the Radon-Nikodym derivative of the Dirac measure $\delta_t$, i.e., ${\rm d}\delta_t= \delta(x-t){\rm d}x$ customarily.

If $\tau_q y$ in $({\bf E}_q)$ is in the limit-point case at $+\infty$, we consider a system interacting with a $\delta$-function potential, given by
\begin{equation}\label{eqn:perturb-probl-lim-pt}
\tau^{t,r}_q y:=-y''(x)+[q(x)-r\delta(x-t)]y(x)=\lambda y(x),\ x\in[0,+\infty),\ y(0)=0,
\end{equation}
while, if $\tau_q y$ in $({\bf E}_q)$ is in the limit-circle case at $+\infty$, the $\delta$-point interaction system, similar to \eqref{eqn:perturb-probl-lim-pt},  is given by
\begin{equation}\label{eqn:perturb-probl-lim-cir}
\begin{split}
&\tau^{t,r}_q y:=-y''(x)+[q(x)-r\delta(x-t)]y(x)=\lambda y(x),\ x\in[0,+\infty),  \\
&y(0)=0=[y,f](+\infty)\sin\beta-[y,g](+\infty)\cos\beta,
\end{split}
\end{equation}
where $q$ is just the one in \eqref{eqn:main-probl}, the pair $(t,r)\in [0,+\infty)\times [0,\epsilon)$, the positive number $\epsilon$ is sufficiently small, and $\beta,\ f,\ g$ are all defined as in \eqref{def-oper-lim-circle}.
The non-negative number $r$ above is called {\bf coupling constant} that stands for the intensity of the interaction.

As  $t$ (respectively, $r$) runs over all real numbers in $[0,+\infty)$ (respectively, $[0, \epsilon]$),  we obtain two families of {\bf perturbation problems} $({\bf E}^{t,r}_{lp})$  and  $({\bf E}^{t,r}_{lc})$, denoted by \eqref{eqn:perturb-probl-lim-pt} and \eqref{eqn:perturb-probl-lim-cir} respectively, of the original problem $({\bf E}_q)$.
The domain of the corresponding differential operators $S^{t,r}_{lp}$ and $S^{t,r}_{lc}$ are respectively given by
\begin{equation}\label{def-pert-oper}
\mathcal{D}(S^{t,r}_{lp}):=\left\{y\in L^{2}[0,+\infty):\ y\in AC_{loc}([0,t)\cup(t,+\infty)),\ \exists y'(t\pm0),\ y(0)=0,\atop
 y'(t-0)-y'(t+0)=ry(t),\ \tau^{t,r}_q y\in L^{2}[0,+\infty)\right\},
\end{equation}
and
\begin{equation}\label{def-pert-oper-lim-cir}
\mathcal{D}(S^{t,r}_{lc}):=\left\{y\in L^{2}[0,+\infty):\ y\in AC_{loc}([0,t)\cup(t,+\infty)),\ \exists y'(t\pm0),\ y'(t-0)-y'(t+0)=ry(t),\atop
y(0)=0=[y,f](+\infty)\sin\beta-[y,g](+\infty)\cos\beta,\ \tau^{t,r}_q y\in L^{2}[0,+\infty)\right\}.
\end{equation}
See \cite{swz2007} and \cite{et2013}.

The next step is to show that, under the hypothesis ${\bf (H_1)}$ (respectively, ${\bf (H_2)}$) proposed above, the perturbation operators $S^{t,r}_{lp}$ (respectively, $S^{t,r}_{lc}$) also have the first eigenvalues $\lambda(t,r;q)$ that are also principle (see Lemma~\ref{lem:dis-spc-Str-lp}, Lemma~\ref{lem:dis-spc-Str-lc} and Theorem~\ref{thm:property-1st-eigenv-funct}). And then suppose that $\lambda(t,r;q)$ are known for any $(t, r)\in [0,+\infty)\times [0,\epsilon]$, which means that the first energy eigenvalues of all perturbation systems can be observed, where $\epsilon$ is sufficiently small. We will show that $\lambda(t,r;q)$, which called {\bf the first eigenvalue function} of $({\bf E}_q)$ as in Definition \ref{defn:1st-eigenv-funct}, is continuous and differentiable with respect to $(t, r)\in[0,+\infty)\times[0,\epsilon]$.

At last,  as Theorem~\ref{thm:main-thm} says, the desired potential $q$ of the main problem $({\bf E}_q)$ can be uniquely and directly reconstructed by
\begin{equation}\label{eqn: q-recover}
q(x)=\frac{\varphi''_0(x)}{\varphi_0(x)}+\lambda_1,
\end{equation}
where $\lambda_1=\lambda(x,0,q)$ and $\varphi_0(x)=\sqrt{-\frac{\partial \lambda(x,0;q)}{\partial r}}$ on $(0,+\infty)$. Note that $\lambda(x,0;q)$ always equals to the first eigenvalue $\lambda_1(q)$ of $({\bf E}_q)$, since in this case each of the problems $({\bf E}^{t,0}_{lp})$ and $({\bf E}^{t,0}_{lc})$ is actually $({\bf E}_q)$.

\bigskip
This paper is organized as follows.
In Section~\ref{sec:prelim}, as a preliminary, we will present a list of the conditions that make sure the existence of the first eigenvalue of $({\bf \tilde{E}}_q)$, and introduce some basic results about the eigenvalues of perturbation problems. These will be used to prove our main theorem.
In Section~\ref{sec:1st-eigen-funct}, we will define and study the first eigenvalue function in the singular case.
Section~\ref{sec:potential-recovery} is the core part, wherein we give the unique reconstruction result of potential by the first eigenvalue function, and some further discussions are also given in the remarks.

\bigskip

\section{Notations and preliminary}\label{sec:prelim}
\medskip

In this section, we firstly give a set of the conditions which ensure that $({\bf \tilde{E}}_q)$ has the principle eigenvalue.
Note that, if all the points in the spectrum of an operator are isolated, we speak of a {\bf pure point spectrum} in the present paper.

\begin{lem}\label{lem:principle-eigenv-exist-case}
$(i)$ If $ q(x)\in L^1_{loc}[0,+\infty)$ is bounded below, $q_0:=\inf_{x\in[0,+\infty)} q(x) > -\infty$, then
$\tau_q y$ in $({\bf \tilde{E}}_q)$ is in the limit-point case at $+\infty$ (see \cite[Theorem 10.1.4]{Hille1969} or \cite[Theorem 7.4.1]{Zettl2005}) and there exists essential spectrum involving in $[q_0,+\infty)$ (see \cite[Theorem 10.3.3]{Hille1969}). Especially if $\alpha=0$ or $\alpha\in[\frac{\pi}{2},\pi)$, there exists a eigenvalue below $q_0$ (see \cite[Theorem 10.3.3]{Hille1969}), which is also the principle eigenvalue (see \cite[Theorem 10.12.1(8)(iii)]{Zettl2005});

\mednoind
$(ii)$ if $q(x)\in L^1_{loc}[0,+\infty)$ is bounded below, and also tends to $+\infty$ as $x$ goes to $+\infty$, then it is obviously seen from the statement $(i)$ that $\tau_q y$ in $({\bf \tilde{E}}_q)$ is in the limit-point case at $+\infty$, and meanwhile the operator $\tilde{S}_q$ has pure point spectrum which is bounded below (see \cite[Theorem 10.3.4]{Hille1969}).
More precisely, in this case, the spectral set of $\tilde{S}_q$ only has isolated real eigenvalues. Let $\lambda_n(q)$ be the $n$-th eigenvalue of $\tilde{S}_q$ or $({\bf \tilde{E}}_q)$, then there holds $-\infty<\lambda_1(q)<\lambda_2(q)<\cdots <\lambda_n(q)\rightarrow +\infty$ as $n\rightarrow +\infty$, and the corresponding eigenfunction associated to $\lambda_n(q)$ has exactly $n-1$ zeros on $(0, +\infty)$ (see \cite[Theorem 10.12.1(8)(ii)]{Zettl2005});

\mednoind
$(iii)$ if $\tau_q y$ in $({\bf \tilde{E}}_q)$ is in the limit-circle case at $+\infty$, then all of the spectral points of $({\bf \tilde{E}}_q)$ are automatically isolated (see \cite[Theorem 10.12.1(2)]{Zettl2005}). In addition, if $\tau_q y$ is  non-oscillatory at $+\infty$, then the spectrum is also bounded below (see \cite[Definition 6.2.2 and Theorem 10.12.1(3)]{Zettl2005}), and, under the separation boundary condition, the principle eigenvalue is exactly the first eigenvalue (see \cite[Theorem 10.12.1(4)]{Zettl2005}).
\end{lem}

\bigskip

In the rest of this section, we will show the existence of the eigenvalue associated to perturbation problems.  For this, we firstly introduce some concepts, terminologies, and several useful facts. For more details, the reader may refer to \cite[Chapter VI]{Kato1980} and \cite[Section 2]{SS1999}, etc.

Let $\BC$ (respectively, $\BR$) be the field of complex (respectively, real) numbers, and $\mathcal{D}$ be a subspace of a Hilbert space $\KH$ over $\BC$.
A mapping $\form{t}[u,v]: \mathcal{D} \times \mathcal{D}\rightarrow \mathbb{C}$ is
called a \emph{sesquilinear form} on $\KH$ if it is linear in $u\in \mathcal{D}$ and semilinear in $v\in \mathcal{D}$. Here $\mathcal{D}$ will be called the \emph{domain} of $\form{t}$, denoted by $\mathcal{D}(\form{t})$, and $\form{t}[u]:=\form{t}[u,u]$ will be called the \emph{quadratic form} associated with $\form{t}[u,v]$. We shall call $\form{t}[u,v]$ or $\form{t}[u]$  simply a {\bf form} when there is no possibility of confusion in the following paragraphs (see \cite[Page 308]{Kato1980}).

A form $\form{t}$ is said to be \emph{symmetric} if
$$\form{t}[u,v]=\overline{\form{t}[v,u]},\ u,v\in \mathcal{D}(\form{t}).$$
A symmetric form $\form{t}$ is said to be \emph{bounded from below} if
$$\form{t}[u]\geq \gamma \|u\|^{2},\ u\in \mathcal{D}(\form{t}),$$
where $\gamma \in \mathbb{R}$ is a constant (see \cite[Page 309]{Kato1980}).
For symmetric and bounded from below form $\form{t}$, a form $\form{s}$ in $\KH$ is said to be \emph{relatively bounded with respect to $\form{t}$}, or simply $\form{t}$-bounded, if $\mathcal{D}(\form{t})\subset\mathcal{D}(\form{s})$ and
$$|\form{s}[u]|\leq a \|u\|^2+b |\form{t}[u]|,\ u\in \mathcal{D}(\form{t}),$$
where $a,\ b$ are non-negative constants. The greatest lower bound for all
possible values of $b$ will be called the \emph{$\form{t}$-bound of $\form{s}$} (see \cite[Page 319]{Kato1980}).

Following the similar arguments in \cite[Page 344]{Kato1980}, we can define two symmetric forms
\begin{equation}\label{form}
 \form{t}^{lp}[u,v]=\int^{+\infty}_{0}(u'\overline{v}'+qu\overline{v})\,\mathrm{d} x,
\end{equation}
and
\begin{equation}\label{forms}
\form{s}^{lp}[u,v]=\int^{+\infty}_{0}(u'\overline{v}'+qu\overline{v})\,\mathrm{d} x-ru(t)\overline{v}(t),
\end{equation}
where $0\leq q(x)\in L^1_{loc}[0,+\infty)$ and
$$u,v\in \mathcal{D}(\form{t}^{lp})=\mathcal{D}(\form{s}^{lp}):=
\left\{y\in L^{2}[0,+\infty):\ y\in AC[0,+\infty),\ y'\in L^{2}[0,+\infty),\ \atop \int_0^{+\infty}q|y|^2\,{\rm d}x<+\infty,\  y(0)=0\right\}.$$
Since $q$ is non-negative, the solution $y$ of $({\bf E}_q)$ satisfies
$y'\in L^2[0,+\infty)$ as well as $|q|^{\frac{1}{2}}y\in L^2[0,+\infty)$ (see \cite[Theorem 10.1.4]{Hille1969}). Therefore, according to the representation theorem (see \cite[Chapter VI, Theorem 2.1]{Kato1980} and \cite[Chapter VI, Theorem 4.2]{Kato1980}), we can know that $S_q$ (respectively, $S^{t,r}_{lp}$) defined in \eqref{def-oper-lim-point} (respectively, \eqref{def-pert-oper}) is in fact the self-adjoint operator associated with the form $\form{t}^{lp}$ (respectively, $\form{s}^{lp}$).

If $\tau_q y$ in $({\bf E}_q)$ is in the limit-circle case and non-oscillatory at $+\infty$, we will give the corresponding forms in the following.

Set $\mathcal{G}[0,+\infty):=\left\{ y\in AC_{loc}[0,+\infty)\, :\, y'\in AC_{loc}[0,+\infty) \right\}.$

Firstly, we can claim that, there exists a strictly positive function $g\in\mathcal{G}[0,+\infty)$ such that the problems $({\bf E}_q)$ and $({\bf E}^{t,r}_{lc})$ are respectively equivalent to the following problems
\begin{equation}\label{eqn:lc}
\begin{split}
&\nu z:=-(Pz')'(x)+Q(x)z(x)=\lambda W(x)z(x),\ x\in[0,+\infty), \\ &z(0)=0=z(+\infty)\sin\beta-Pz'(+\infty)\cos\beta,
\end{split}
\end{equation}
and
\begin{equation}\label{eqn:perturb-lc}
\begin{split}
&\nu^{t,r}z:= -(Pz')'(x)+[Q(x)-r\delta(x-t)g^2(x)]z(x)=\lambda W(x)z(x),\ x\in[0,+\infty),  \\
&y(0)=0=z(+\infty)\sin\beta-Pz'(+\infty)\cos\beta,
\end{split}
\end{equation}
where $z=\frac{y}{g}$, $P=W=g^2,\ Q=g\cdot \tau_q g$ and $\beta$ is defined as in \eqref{def-oper-lim-circle}.

Indeed, since $\tau_q y$ in $({\bf E}_q)$ is in the limit-circle case and non-oscillatory at $+\infty$, by \cite[Theorem 8.2.1]{Zettl2005},
there exists a real-valued principal solution $u_{+\infty}$ and a real-valued non-principal solution $v_{+\infty}$ of $\tau_q y=\lambda_0 y$ at $+\infty$ for some real $\lambda_0$ satisfying $v_{+\infty}>0$ and $[u_{+\infty},\ v_{+\infty}](+\infty)=1$.
Here the word ``{\bf principal}" about $u_{+\infty}$ at $+\infty$ means that, there exists some $a\in(0,+\infty)$ such that $u_{+\infty}(x)\neq 0$ for $x\in[a,+\infty)$, and every solution $y$ of $\tau_q y=\lambda_0 y$ which is not a multiple of $u_{+\infty}$ satisfies $\lim\limits_{x\rightarrow +\infty}\frac{u_{+\infty}(x)}{y(x)}=0$, while the word ``{\bf non-principal}" about $v_{+\infty}$ at $+\infty$ means that, $v_{+\infty}$ is not principal at $+\infty$ and is non-zero on $[a,+\infty)$ for some $a\in(0,+\infty)$. Note that there exists some $c\in(0,+\infty)$ such that $u_{+\infty}(x)v_{+\infty}(x)\neq 0$ for $x\in[c,+\infty)$. Hence, we can see that $u_{+\infty}$ and $v_{+\infty}$ are a pair of (BC) solutions defined in \eqref{def-oper-lim-circle}.
Consequently, setting $g=v_{+\infty}$ and $z=\frac{y}{g}$, by \cite[Lemma 3.2, Lemma 3.7(3) and Lemma 4.3]{nz1992} (also see \cite[Theorem 8.3.1]{Zettl2005}),
it can be directly seen that the problems $({\bf E}_q)$ and $({\bf E}^{t,r}_{lc})$ are respectively equivalent to \eqref{eqn:lc} and \eqref{eqn:perturb-lc}, namely, the claim above is ture.

The domain of the differential operators $T_{lc}$ and $T^{t,r}_{lc}$ related to \eqref{eqn:lc} and \eqref{eqn:perturb-lc} are respectively given by
\begin{equation}\label{def-lc-oper}
\mathcal{D}(T_{lc}):=\left\{z\in L_{W}^{2}[0,+\infty):\ z,\, Pz'\in AC_{loc}[0,+\infty),\ \nu z\in L_{W}^2[0,+\infty),\atop z(0)=0=z(+\infty)\sin\beta-Pz'(+\infty)\cos\beta \right\},
\end{equation}
and
\begin{equation}\label{def-pert-lc-oper}
\mathcal{D}(T^{t,r}_{lc}):=\left\{z\in L_{W}^{2}[0,+\infty):\ z\in AC_{loc}([0,t)\cup(t,+\infty)),\ \nu^{t,r}z\in L_{W}^{2}[0,+\infty),
\exists Pz'(t\pm 0),\atop Pz'(t-0)-Pz'(t+0)=rg^2(t)z(t),\
z(0)=0=z(+\infty)\sin\beta-Pz'(+\infty)\cos\beta \right\}.
\end{equation}
Here $L_{W}^2[0,+\infty)$ denotes the space of real-valued weighted square-integrable functions on $[0,+\infty)$ equipped with the canonical norm $\norm{\cdot}_W:=\int_0^{+\infty} W|\cdot|^2 {\rm d}x$.

Then, we can define two symmetric forms
\begin{equation}\label{form-lc}
 \form{t}^{lc}[u,v]=\tan \beta uv(+\infty)+\int^{+\infty}_{0}(Pu'\overline{v}'+Qu\overline{v})\,\mathrm{d} x,
\end{equation}
and
\begin{equation}\label{forms-lc}
\form{s}^{lc}[u,v]=\tan \beta uv(+\infty)+\int^{+\infty}_{0}(Pu'\overline{v}'+Qu\overline{v})\,\mathrm{d} x-r g^2(t)u(t)\overline{v}(t),
\end{equation}
where
$$u,v\in \mathcal{D}(\form{t}^{lc})=\mathcal{D}(\form{s}^{lc}):=
\left\{z\in L_{W}^{2}[0,+\infty):\ z\in AC[0,+\infty),\ \sqrt{P}z'\in L^{2}[0,+\infty),\ \atop \int_0^{+\infty}Q|z|^2\,{\rm d}x<+\infty,\  z(0)=0\right\}.$$
Since $\sqrt{P}z'\in L^{2}[0,+\infty)$ (see \cite[Lemma 3.7]{Zettl2005}) and $z$ is bounded,  according to the representation theorem (see \cite[Chapter VI, Theorem 2.1]{Kato1980}), we can know that $T_{lc}$ (respectively, $T^{t,r}_{lc}$) is in fact the self-adjoint operator associated with the form $\form{t}^{lc}$ (respectively, $\form{s}^{lc}$). In other words, the forms $\form{t}^{lc}$  and $\form{s}^{lc}$ are the forms associated to the problems $({\bf E}_q)$ and $({\bf E}^{t,r}_{lc})$, respectively.

\bigskip

Now we list some notations introduced in \cite[Page 202]{Kato1980}. Consider two closed linear manifolds $M$ and $N $ of a Banach space $X$. Denote by $S_M$ the unit sphere of $M$. Set
$$
\delta(M, N) := \sup_{u \in S_M} \text{dist}(u, N), \quad \hat{\delta}(M, N) := \max[\delta(M, N), \delta(N, M)].
$$
Consider the set $\mathscr{C}(X, Y)$ of all closed operators from a Banach space $X $ to a Banach space $ Y $. If $ T, S \in \mathscr{C}(X, Y) $, then their graphs $ G(T), G(S) $ are closed linear manifolds of the product space $ X \times Y$. Set
$$
\delta(T, S) := \delta(G(T), G(S)), \quad \hat{\delta}(T, S) := \hat{\delta}(G(T), G(S)) = \max[\delta(T, S), \delta(S, T)].
$$
$\hat{\delta}(T, S)$ will be called the \emph{gap between $T$ and $S$}. We recall that $T_n $ converges to $ T$ \textbf{in the generalized sense} if $\hat{\delta}(T_n, T) \rightarrow 0$.

Then we recall the definition of the norm resolvent convergence.
Denote by $\BN$ the set of non-negative integers.
Let $T$ and $\{T_m\}_{m\in\BN}$ be closed operators on a Hilbert space over $\BC$. The sequence $\{T_m\}_{m\in\BN}$ converges to $T$ in the sense of {\bf norm resolvent convergence}, denoted by $T_m \stackrel{R}{\Rightarrow} T,$
if there is a number $\mu \in \BC$ belonging to the resolvent sets $\rho(T)$ and $\rho(T_m)$ for all $m\in \BN$ and the sequence of bounded operators $(T_m-\mu)^{-1}$ converge uniformly to the operator $(T-\mu)^{-1}$ as $m\rightarrow +\infty$ (see \cite[Section 2]{SS1999} or \cite[VIII, Page 427]{Kato1980}).

\begin{lem}\label{lem: converg-eigenv}
For self-adjoint operators $\{T_m\}_{m\in\BN}$ and $T$, suppose that $T_m \stackrel{R}{\Rightarrow} T,$\\
$(i)$ if $\lambda(m)\in\sigma(T_m)$, and $\lambda(m)\rightarrow c$ as $m\rightarrow +\infty,$ then $c\in\sigma(T)$ (see \cite[Corollary 1.4]{ysz2017});\\
$(ii)$ if $\lambda\in\sigma(T)$, then there must exist $\lambda(m)\in\sigma(T_m),$ such that $\lambda(m)\rightarrow \lambda$ as $m\rightarrow +\infty$ (see \cite[Corollary 1.4]{ysz2017}). If $\lambda$  is an isolated eigenvalue of the operator $T$ with finite geometric multiplicity $n$, then there are finitely many eigenvalues of the operators $T_m$ in an arbitrary sufficiently small $\delta$ neighborhood of the point $\lambda$ if $m$ is large enough. Moreover, their total geometric multiplicity equals $n$ (see \cite[Lemma 1.5]{ysz2017});\\
$(iii)$ the $n$-th eigenvalue of $T$ is continuous under perturbations in the sense of the norm resolvent convergence, that is, $\lambda_n^m\rightarrow \lambda_n$ as $m\rightarrow +\infty$, where $\lambda_n^m$ and $\lambda_n$ are respectively the $n$-th eigenvalues of $T_m$ and $T$ (see \cite[Lemma~1]{YSh2014}, \cite[Lemma 1.6]{ysz2017} or \cite[Theorem 6]{SS1999}).
\end{lem}

\begin{lem}\label{lem:dis-spc-Str-lp}
For the problem $({\bf E}_q)$, if the potential $q$ satisfies the condition in Lemma~\ref{lem:principle-eigenv-exist-case}$(i)$, then the self-adjoint operator $S^{t,r}_{lp}$ has eigenvalues that are bounded below.
In particular, if this $q$ satisfies the condition in Lemma~\ref{lem:principle-eigenv-exist-case}$(ii)$, the perturbed operator $S^{t,r}_{lp}$ has pure point spectrum, namely, its spectrum has only isolated real eigenvalues which can be ordered as in Lemma~\ref{lem:principle-eigenv-exist-case}$(ii)$.
\end{lem}
\begin{prf}  If $q$ in $({\bf E}_q)$ belongs to the case $(i)$ in Lemma \ref{lem:principle-eigenv-exist-case}, then, without loss of generality, we may as well assume that $q$ is non-negative, so the essential spectra of $({\bf E}^{t,r}_{lp})$ is the same as $({\bf E}_q)$ (see \cite[Theorem 2.26]{km2013} or \cite[Theorem 4]{ys2014}).

For any given $(t,r)\in [0,+\infty)\times [0,+\infty)$, consider $T:=S_q$ and $T_m:=S^{t_m,r_m}_{lp}$, where $\{(t_m,r_m)\}_{m\in\BN}\subset [0,+\infty)\times [0,+\infty)$ tends to $(0,0)$.
And let $\form{t}^{lp}$ and $\form{s}_m^{lp}$ be the forms related to $T$ and $T_m$, defined as \eqref{form} and \eqref{forms}, respectively.
According to \cite[Chapter IV, (1.19)]{Kato1980}, for any $\epsilon>0,$ there exists $\Gamma(\epsilon)>0$ such that
\begin{equation}\label{cauchy}
\|y\|_{\infty}^{2}\leq
\epsilon\|y'\|_{2}^{2}+\Gamma(\epsilon)\|y\|_{2}^{2},\ y\in\left\{y\in L^{2}[0,+\infty):\ y\in AC[0,+\infty),\ y'\in L^{2}[0,+\infty) \right\}.
\end{equation}
Note that $\mathcal{D}(\form{t}^{lp})=\mathcal{D}(\form{s}_m^{lp})$, and hence for any $u\in \mathcal{D}(\form{t}^{lp}),\ \epsilon>0$, there exists $\Gamma(\epsilon)>0$ such that
$$|(\form{t}^{lp}-\form{s}_m^{lp})[u]|=r_m |u|^2(t_m)\leq r_m\epsilon\|u'\|_{2}^{2}+r_m\Gamma(\epsilon)\|u\|_{2}^{2}\atop
\leq
r_m\epsilon(\int^{+\infty}_{0} \left(|u'|^{2}+q|u|^{2}\right)\mathrm{d}x)+r_m\Gamma(\epsilon)\|u\|_{2}^{2}
=r_m\epsilon \form{t}^{lp}[u]+r_m\Gamma(\epsilon)\|u\|_{2}^{2}.$$
Since $r_m\rightarrow 0$ as $m\rightarrow +\infty$, by \cite[Chapter VI, Theorem 3.6]{Kato1980}, we know that the operators $T_m$ converges to $T$  in the generalized sense. Then, according to \cite[Chapter IV, Theorem 2.23(b)(c)]{Kato1980}, we can see that $T_m$ converges to $T$  in the generalized sense is equivalent to $T_m \stackrel{R}{\Rightarrow} T$, and then, by Lemma \ref{lem: converg-eigenv}$(ii)$ and Lemma \ref{lem:principle-eigenv-exist-case}$(i)$, for any $(t,r)\in [0,+\infty)\times [0,+\infty)$, it follows that $S^{t,r}_{lp}$ has an eigenvalue.

Moreover, since $q$ is bounded below, the spectrum of $S^{t,r}_{lp}$ is bounded below.
Indeed, according to \eqref{cauchy} again, for sufficiently small $\epsilon$, one can verify that
\begin{equation}\label{bounded-below}
\form{s}^{lp}[u]=\form{t}^{lp}[u]+\form{h}^{lp}[u]
 =\int^{+\infty}_{0} \left(|u'|^{2}+q|u|^{2}\right)\mathrm{d}x-r|u(t)|^2\atop
 \geq(1-r\epsilon)\|u'\|^2-r\Gamma(\epsilon)\|u\|^2\geq -r\Gamma(\epsilon)\|u\|^2,
\end{equation}
and hence the spectrum of $S^{t,r}_{lp}$ is bounded below.

Now, we consider that $q$ satisfies the condition in Lemma~\ref{lem:principle-eigenv-exist-case}$(ii)$.
It is well known that the self-adjoint operator has  pure point spectrum if and only if its resolvent operator is compact (see \cite[Chapter IX, Theorem 3.1]{EE1987}).
Therefore, to prove $S^{t,r}_{lp}$ has pure point spectrum, it just needs to prove the resolvent operator of $S^{t,r}_{lp}$ is compact.

For fixed $t\in[0,+\infty)$, set a form
$$\form{h}^{lp}[u,v]=-ru(t)\overline{v}(t),$$
where $u,v\in \mathcal{D}(\form{h}^{lp}):=
\left\{y\in L^{2}[0,+\infty):\ y\in AC[0,+\infty)\right\}$.
To show $S^{t,r}_q$ has compact resolvent, by \cite[Chapter VI, Theorem 3.4]{Kato1980}, we need to prove that $\form{h}^{lp}$ is $\form{t}^{lp}$-bounded form, and $\form{t}^{lp}$-bound is less than $\frac{1}{2}$.

According to \eqref{cauchy}, for any $b>0$, there exists $a>0$ such that
\begin{align*}
b\form{t}^{lp}[u]-\form{h}^{lp}[u]+a\|u\|^2
& =b\int^{+\infty}_{0} \left(|u'|^{2}+q|u|^{2}\right)\mathrm{d}x+r|u(t)|^2+a\|u\|^2\\
& \geq(b+r\epsilon)\|u'\|^2+(a+r\Gamma(\epsilon))\|u\|^2\geq 0,
\end{align*}
which means $\form{t}^{lp}-$bound of $\form{h}^{lp}$ is 0.
Hence, we know that $S^{t,r}_{lp}$ has pure point spectrum.
\end{prf}

\begin{lem}\label{lem:dis-spc-Str-lc}
If $\tau_q y$ in $({\bf E}_q)$ is in the limit-circle case and non-oscillatory at $+\infty$, the self-adjoint operator $S^{t,r}_{lc}$ has pure point spectrum that is bounded below, and all eigenvalues are simple.
\end{lem}

\begin{prf} It directly follows from Lemma \ref{lem:principle-eigenv-exist-case}$(iii)$, \cite[Theorem 6]{zhang2014} and \cite[Corollary 8.4]{et2013}.
\end{prf}

\bigskip
\section{The first eigenvalue function $\lambda(t,r)$}\label{sec:1st-eigen-funct}
\medskip

For the perturbation problems $({\bf E}_{lp}^{t,r})$ and $({\bf E}_{lc}^{t,r})$ with $(r,t)\in [0, +\infty)\times [0, +\infty)$, we suppose that the potential $q\in L^1_{loc}[0,+\infty)$ is unknown, while the first eigenvalue is known.
By moving the position $t$ of interaction along the open interval $(0,+\infty)$ and adjusting the intensity $r$, we gain a family of the first eigenvalues, denoted by $\lambda(t,r;q)$, of $({\bf E}_{lp}^{t,r})$ or $({\bf E}_{lc}^{t,r})$ for all $(r,t)\in (0, +\infty)\times (0, +\infty)$.

\begin{defn}\label{defn:1st-eigenv-funct}
The first eigenvalues $\lambda(t,r;q)$ above is called {\bf the first eigenvalue function} of the original problem $({\bf E}_q)$.
\end{defn}

\noindent
{\bf Notation:} For the given problem $({\bf E}_q)$, the potential $q$ has already been fixed, although it is unknown. Hence, {\bf for convenience, thereinafter we sometimes simply denote by $\lambda(t,r)$ the first eigenvalue function when there is no risk of confusion}.
\medskip

In the packet of theorems below (i.e., Theorems~\ref{thm:cont-1st-eigenv-funct}, \ref{thm:property-1st-eigenv-funct} and \ref{thm:diff-1st-eigenv-funct}), we will prove some essential properties of $\lambda(t,r;q)$, such as continuity, boundedness and differentiability, etc.
These properties for the regular S-L problems had been presented in \cite{ZQC2024} and the references therein.
Maybe some properties can be acquired analogously from some references, but we have not found the related evidences explicitly stated in the literature. Hence, for the completeness and rigorousness of our paper as well as the convenience for readers, we will give self-contained and explicit proofs of these facts.

\begin{thm}\label{thm:cont-1st-eigenv-funct}
If $({\bf E}_q)$ satisfies either ${\bf (H_1)}$ or ${\bf (H_2)}$, then $\lambda(t,r)$ is continuous with $(t,r)\in [0,+\infty)\times [0,+\infty).$
\end{thm}

\begin{prf}
For any given $(t,r)\in [0,+\infty)\times [0,+\infty)$, consider $T_m:=S^{t_m,r_m}_{lp}$ and $T:=S^{t,r}_{lp}$, or $T_m:=S^{t_m,r_m}_{lc}$ and $T:=S^{t,r}_{lc}$, where $\{(t_m,r_m)\}_{m\in\BN}\subset [0,+\infty)\times [0,+\infty)$ tends to $(t,r)$.

If $({\bf E}_q)$ satisfies ${\bf (H_1)}$,
let $\form{s}_{t_m,r_m}^{lp}$ and $\form{s}_{t,r}^{lp}$ respectively be the forms related to $T_m:=S^{t_m,r_m}_{lp}$ and $T:=S^{t,r}_{lp}$, defined as \eqref{forms}.
For $u\in \mathcal{D}(\form{s}_{t_m,r_m}^{lp})=\mathcal{D}(\form{s}_{t,r}^{lp})$, by \eqref{cauchy}, there exists $\epsilon,\Gamma(\epsilon)>0$ such that
\begin{equation}\label{equ:r-bdd-lp}
|(\form{s}_{t_m,r_m}^{lp}-\form{s}_{t,r}^{lp})[u]|=r_m |u|^2(t_m)-r |u|^2(t)\leq (r_m-r)(\epsilon\|u'\|_{2}^{2}+\Gamma(\epsilon)\|u\|_{2}^{2})\rightarrow 0,\ m \rightarrow +\infty.
\end{equation}
It follows from the Max-Min principle of forms (see \cite[Theorem XIII.2]{RS1978}) that
\begin{equation*}
\lambda(t,r)=\inf\left\{\frac{\form{s}_{t,r}^{lp}[u]}{\|u\|^2}:\ u\in \mathcal{D}(\form{s}_{t,r}^{lp})\right\},\ \lambda(t_m,r_m)=\inf\left\{\frac{\form{s}_{t_m,r_m}^{lp}[u]}{\|u\|^2}:\ u\in \mathcal{D}(\form{s}_{t_m,r_m}^{lp})\right\},
\end{equation*}
and the continuity of $\lambda(t,r)$ of $S^{t,r}_{lp}$ holds.

If $({\bf E}_q)$ satisfies ${\bf (H_2)}$,
let $\form{s}_{t_m,r_m}^{lc}$ and $\form{s}_{t,r}^{lc}$ respectively be the forms related to $T_m:=S^{t_m,r_m}_{lc}$ and $T:=S^{t,r}_{lc}$, defined as \eqref{forms-lc}.
For any $u\in \mathcal{D}(\form{s}_{t_m,r_m}^{lc})=\mathcal{D}(\form{s}_{t,r}^{lc})$, it is apparent that $u\in L_{W}^2$, where $W=g^2$. Then, there holds
\begin{equation}\label{equ:r-bdd-lc}
|(\form{s}_{t_m,r_m}^{lc}-\form{s}_{t,r}^{lc})[u]|=r_m |gu|^2(t_m)-r |gu|^2(t)\leq (r_m-r)\|gu\|_{\infty}^{2}\rightarrow 0,\ m \rightarrow +\infty.
\end{equation}
It follows from the Max-Min principle of forms (see \cite[Theorem XIII.2]{RS1978}) that
\begin{equation*}
\lambda(t,r)=\inf\left\{\frac{\form{s}_{t,r}^{lc}[u]}{\|u\|_W^2}:\ u\in \mathcal{D}(\form{s}_{t,r}^{lc})\right\},\ \lambda(t_m,r_m)=\inf\left\{\frac{\form{s}_{t_m,r_m}^{lc}[u]}{\|u\|_W^2}:\ u\in \mathcal{D}(\form{s}_{t_m,r_m}^{lc})\right\},
\end{equation*}
and the continuity of $\lambda(t,r)$ of $S^{t,r}_{lc}$ also holds.
\end{prf}

\begin{rem}\label{rem:bdd-r-lambda}
From \eqref{equ:r-bdd-lp} and \eqref{equ:r-bdd-lc} in the proof of Theorem \ref{thm:cont-1st-eigenv-funct}, we can see that if $r$ is bounded, then $\lambda(t,r)$ is bounded.
\end{rem}

For further studying the function $\lambda(t,r)$, we recall the concept of {\bf Weyl function}.
Let $\varphi(x,\lambda)$ and $\psi(x,\lambda)$ be solutions
of $-y''+qy=\lambda y$ such that
\begin{equation}\label{315}
 \varphi(0)=0, \ \varphi'(0)=1;\  \psi(0)=1, \ \psi'(0)=0.
\end{equation}
Denote by $m(\lambda)$  the Weyl function of differential expression $\tau_q y$. All conclusions in the next lemma follow from the classical Weyl theory. For more details, the reader may refer to \cite{Hille1969}.
\begin{lem}\label{weyl}
$(i)$ For any $\lambda$ in the resolvent set $\rho(S_q)$ of $S_q$, if $\tau_q y$ in $({\bf E}_q)$ is in the limit-point case at $+\infty$, there exists $m(\lambda)$ such that $\chi(x,\lambda)=\psi(x,\lambda)+m(\lambda)\varphi(x,\lambda)\in L^2(0,+\infty)$ (see \cite[Page 498]{Hille1969}).
If $\tau_q y$ in $({\bf E}_q)$ is in the limit-circle case at $+\infty$,  there exists $m(\lambda)$ such that $\chi(x,\lambda)=\psi(x,\lambda)+m(\lambda)\varphi(x,\lambda)\in L^2(0,+\infty)$ satisfying $[\chi,f](+\infty)\sin\beta-[\chi,g](+\infty)\cos\beta=0$.\\
$(ii)$ $m'(\lambda)=\int^{+\infty}_0 \chi^2(x,\lambda) \mathrm{d}x$ (see \cite[Page 523(10.3.14)]{Hille1969}). \\
$(iii)$ $m(\lambda)$ and $\chi(x,\lambda)$ are both holomorphic with respect to $\lambda$ (see \cite[Lemma 10.3.1]{Hille1969}), and $m(\lambda)$ only has simple zero point and pole on $\BR$ (see \cite[Page 512, Theorem 10.3.1]{Hille1969}). \\
$(iv)$ $\lambda$ is an eigenvalue of the problem $({\bf E}_q)$ if and only if it is the pole of $m(\lambda)$. And $m(\lambda)$ is strictly increasing on the open interval between $-\infty$ and the first eigenvalue and every open interval between any two adjacent eigenvalues (see \cite[Theorem 10.3.1]{Hille1969}).
\end{lem}

\begin{thm}\label{thm:property-1st-eigenv-funct}
Let $\lambda(t,r)$ be the first eigenvalue function of $({\bf E}_q)$.
If $({\bf E}_q)$ satisfies either ${\bf (H_1)}$ or ${\bf (H_2)}$, then it holds that
$$
\lambda(t,r)<\lambda(0,r)=\lambda(t,0)=\lambda_1,\quad (t,r)\in(0,+\infty)\times(0,+\infty),
$$
where $\lambda_1$ is the first eigenvalue of $({\bf E}_q)$.
Meanwhile, the corresponding eigenfunction  $\Phi(x,\lambda(t,r))$ does not change its sign on $(0,+\infty)$ for $(t,r)\in [0,+\infty)\times [0,\epsilon_0]$, provided that $\epsilon_0>0$ is  sufficient small. This means that $\lambda(t,r)$ gives the principle eigenvalues of the perturbation problems of $({\bf E}_q)$.
\end{thm}
\begin{prf} Clearly, $\lambda(t,0)=\lambda_1$ for all $t\in[0,+\infty)$.
If $t=0$, the perturbation term $-r\delta(x-t)$
has no influence to the original problem  $({\bf E}_q)$, and hence $\lambda(0,r)=\lambda_1$.

We firstly prove $\lambda(t,r)\leq\lambda_1$.

Suppose that  $({\bf E}_q)$ satisfies ${\bf (H_1)}$.
Note that $\form{s}^{lp}[u]<\form{t}^{lp}[u]$ for any $u\in \mathcal{D}(\form{t}^{lp})=\mathcal{D}(\form{s}^{lp})$, where $\form{t}^{lp}$ and $\form{s}^{lp}$ are two forms defined as \eqref{form} and \eqref{forms}, respectively. It follows from the Max-Min principle of forms (see \cite[Theorem XIII.2]{RS1978}) that
\begin{equation}\label{max-min}
\lambda(t,r)=\inf\left\{\frac{\form{s}^{lp}[u]}{\|u\|_2^2}:\ u\in \mathcal{D}(\form{s}^{lp})\right\},\ \lambda_1=\inf\left\{\frac{\form{t}^{lp}[u]}{\|u\|_2^2}:\ u\in \mathcal{D}(\form{t}^{lp})\right\},
\end{equation}
and hence $\lambda(t,r)\leq\lambda_1$.

If $({\bf E}_q)$ satisfies ${\bf (H_2)}$, then, for the finite interval $[0,b]$ with $b\in (0,+\infty)$, it follows from the same lines in the proof of \cite[Theorem 3.2]{ZQC2024} that
$\lambda_b(t,r)<\lambda_b$, where $\lambda_b(t,r)$ and $\lambda_b$ are respectively the first eigenvalues of the equation in \eqref{eqn:main-probl} and \eqref{eqn:perturb-probl-lim-cir} with the boundary condition $$y(0)=0=[y,f](b)\sin\beta-[y,g](b)\cos\beta.$$
Furthermore, the singular problem can be approximated by regular problems, that is,  $\lim\limits_{b\rightarrow +\infty}\lambda_b(t,r)=\lambda(t,r)$ (see \cite[Theorem 6]{zhang2014})
and $\lim\limits_{b\rightarrow +\infty}\lambda_b=\lambda_1$ (see \cite[Theorem 10.8.2]{Zettl2005}), and hence $\lambda(t,r)\leq\lambda_1$.

Next, we need to prove that $\lambda(t,r)$ is strictly less than $\lambda_1$ for $(t,r)\in(0,+\infty)\times(0,+\infty)$.

If  $\Phi(x,\lambda(t,r))$ is the first eigenfunction of $({\bf E}^{t,r}_{lp})$, then $\Phi(x,\lambda(t,r))$ satisfies
\begin{equation}\label{eq:lp-1st-eigenf}
-\Phi''+[q-r\delta(x-t)]\Phi=\lambda(t,r) \Phi,\ \Phi(0)=0,
\end{equation}
or equivalently
\begin{equation}\label{eq:lp-1st-eigenf-equiv}
\left\{\aligned
&-\Phi''+q \Phi=\lambda(t,r) \Phi,\ t\neq x\in(0,+\infty),\\
& \Phi'(t-0,\lambda(t,r))-\Phi'(t+0,\lambda(t,r))=r\Phi(t,\lambda(t,r)),\\
& \Phi(0,\lambda(t,r))=0,
\endaligned\right.
\end{equation}
where $\Phi'=\frac{\partial\Phi(x,\lambda(t,r))}{\partial x}$.

If $\Phi(x,\lambda(t,r))$ is the first eigenfunction of $({\bf E}^{t,r}_{lc})$, then $\Phi(x,\lambda(t,r))$ satisfies
\begin{equation}\label{eq:lc-1st-eigenf}
-\Phi''+[q-r\delta(x-t)]\Phi=\lambda(t,r) \Phi,\ \Phi(0)=0=[\Phi,f](+\infty)\sin\beta-[\Phi,g](+\infty)\cos\beta,
\end{equation}
or equivalently
\begin{equation}\label{eq:lc-1st-eigenf-equiv}
\left\{\aligned
&-\Phi''+q \Phi=\lambda(t,r) \Phi,\ t\neq x\in(0,+\infty),\\
& \Phi'(t-0,\lambda(t,r))-\Phi'(t+0,\lambda(t,r))=r\Phi(t,\lambda(t,r)),\\
& \Phi(0,\lambda(t,r))=0=[\Phi,f](+\infty)\sin\beta-[\Phi,g](+\infty)\cos\beta.
\endaligned\right.
\end{equation}
Then, in either case of ${\bf (H_1)}$ and ${\bf (H_2)}$, if $\lambda(t,r)=\lambda_1$, i.e., $\Phi$ is the first eigenfunction of $({\bf E}_q)$, then the second equation in either one of \eqref{eq:lp-1st-eigenf-equiv} and \eqref{eq:lc-1st-eigenf-equiv} can ensure that $\Phi(t)=0$, which contradicts that the first eigenfunction of $({\bf E}_q)$ does not have any zero on $(0,+\infty)$ (see \cite[Theorem 10.12.1(4) and (8)(ii)]{Zettl2005}). And hence $\lambda(t,r)<\lambda_1$ for any $(t,r)\in(0,+\infty)\times(0,+\infty)$.

Finally, we prove that the first eigenfunction $\Phi(x,\lambda(t,r))$  of $({\bf E}^{t,r}_{lp})$ or $({\bf E}^{t,r}_{lc})$ does not change its sign on $(0,+\infty)$.
From Remark \ref{rem:bdd-r-lambda}, we can choose a sufficiently small real number $r$ which is not bigger than $\epsilon_0$, such that
\begin{equation}\label{small}
|\lambda_1-\lambda(t,r)|=|\lambda(t,0)-\lambda(t,r)|
\end{equation}
is small enough for $(t,r)\in[0,+\infty)\times[0,\epsilon_0]$.
According to Lemma \ref{weyl}$(iv)$, it is known that $\lambda_1$ is a pole of $m(\lambda)$, and $m(\lambda)$ is strictly increasing on $(-\infty, \lambda_1)$.
Since  $\lambda(t,r)<\lambda_1$ for any $(t,r)\in(0,+\infty)\times(0,+\infty)$, one can see that
$m(\lambda(t,r))\neq 0$ for any $(t,r)\in(0,+\infty)\times(0,\epsilon_0]$.
Let
\begin{equation}\label{Psi}
\Psi(x,\lambda(t,r))=
\frac{\chi(x,\lambda(t,r))}{m(\lambda(t,r))},\ (t,r)\in(0,+\infty)\times(0,\epsilon_0].
\end{equation}
Then the eigenfunction  $\Phi(x,\lambda(t,r))$ of $({\bf E}^{t,r}_{lp})$ or $({\bf E}^{t,r}_{lc})$ can be chosen as
\begin{equation}\label{eqn:eigenf-perturb-probl}
\Phi(x,\lambda(t,r))=\left\{\aligned
&\varphi(x,\lambda(t,r)), && x\in[0,t],\\
&c(t,r)\Psi(x,\lambda(t,r)),  && x\in[t,+\infty),\\
\endaligned\right.
\end{equation}
where $c(t,r)=\frac{\varphi(t,\lambda(t,r))}{\Psi(t,\lambda(t,r))}$.

Since $\lambda(t,r)<\lambda(t,0)=\lambda_1$ and
the first eigenfunction of $({\bf E}_q)$ associated to $\lambda_1$ does not change sign on $(0,+\infty)$ (see \cite[Theorem 10.12.1(8)(ii)]{Zettl2005}),
it follows from the oscillation theory of linearly second order differential equations (see \cite[Theorem 2.6.2]{Zettl2005})
that
$\varphi(x,\lambda(t,r))$ and $\Psi(x,\lambda(t,r))$ do not change sign on $x\in(0,+\infty)$,  and hence $c(t,r)$ does not change sign, which means that
$\Phi(x,\lambda(t,r))$ does not change its sign on $(0,+\infty)$ for all $(t,r)\in[0,+\infty)\times[0,\epsilon_0]$.
\end{prf}

It is worth noticing that, the statement $\lambda(t,r)<\lambda(t,0)$ in the above theorem, can be deduced by the derivative formula \eqref{eq:part-deriv-r}, but the above proof does not involve this formula.
In other words, the property $\lambda(t,r)<\lambda(t,0)$ is independent of the monotonicity of $\lambda(t,r)$ with respect to $r$.

\begin{lem}\label{thm: equiv-1st-eigenv}
If $({\bf E}_q)$ satisfies either  ${\bf (H_1)}$ or ${\bf (H_2)}$, then $\lambda(t,r)$ is the first eigenvalue function of $({\bf E}_q)$ if and only if
\begin{equation}\label{eqn:equiv-1st-eigenv}
 r\varphi(t,\lambda(t,r))\Psi(t,\lambda(t,r))
 -\frac{1}{m(\lambda(t,r))}=0,\ (t,r)\in(0,+\infty)\times(0,\epsilon_0],
\end{equation}
where $\epsilon_0>0$ is defined as Theorem \ref{thm:property-1st-eigenv-funct} and $\varphi(x,\lambda)$, $\Psi(x,\lambda)$ are given by \eqref{315} and \eqref{Psi}, respectively.
\end{lem}

\noindent
\begin{prf}
For the eigenfunction $\Phi(x,\lambda(t,r))$, defined as \eqref{eqn:eigenf-perturb-probl}, of $({\bf E}^{t,r}_{lp})$ (respectively, $({\bf E}^{t,r}_{lc})$), by the second equality in \eqref{eq:lp-1st-eigenf-equiv} (respectively, \eqref{eq:lc-1st-eigenf-equiv}), we have
$$
\varphi'(t,\lambda)-c\Psi'(t,\lambda)=r\varphi(t,\lambda).
$$
And hence
$$
r\varphi(t,\lambda)=\varphi'(t,\lambda)
-\frac{\varphi(t,\lambda)}{\Psi(t,\lambda)}\Psi'(t,\lambda),
$$
which means \eqref{eqn:equiv-1st-eigenv} holds due to the definition of $\Psi(x,\lambda)$ in \eqref{Psi}. Since every steps in the above argument is reversible, we can know that $\lambda(t,r)$ satisfying \eqref{eqn:equiv-1st-eigenv}
implies $\lambda(t,r)$ being the first eigenvalue function of $({\bf E}_q)$.
\end{prf}

We further obtain the differentiability both of the first eigenvalue function $\lambda(t,r)$ and the corresponding normalized eigenfunction.

\begin{thm}\label{thm:diff-1st-eigenv-funct}
Let $\lambda(t,r)$ be the first eigenvalue function of $({\bf E}_q)$. If $({\bf E}_q)$ satisfies either ${\bf (H_1)}$ or ${\bf (H_2)}$, then we have the following conclusion:
for any $N\in \BN$, there exists a sufficiently small $\epsilon_N\in(0,\epsilon_0]$, such that

\noindent
$(i)$ $\lambda(t,r)$ is continuous differentiable with respect to $(t,r)\in[0,N]\times[0,\epsilon_N]$ and
\begin{equation}\label{equ:abs-cont-1st-eigenv-funct}
\frac{\partial\lambda(t,r)}{\partial t},\ \frac{\partial\lambda(t,r)}{\partial r} \in AC_\BR[0,N],\text{ for any given $r\in [0,\epsilon_N]$}.
\end{equation}
 Moreover,
\begin{equation}\label{eq:part-deriv-r}
\frac{\partial \lambda(t,r)}{\partial r}=-\Phi_0^2(t,\lambda(t,r))\text{ on $[0,+\infty)$},
\end{equation}
where $\Phi_0(x,\lambda(t,r))$ is the normalized eigenfunction of $({\bf E}^{t,r}_{lp})$ or $({\bf E}^{t,r}_{lc})$ associated to $\lambda(t,r)$;

\noindent
$(ii)$ $\Phi(x,\lambda(t,r))$ is continuous differentiable with respect to $r\in [0,\epsilon_N]$.
\end{thm}

\begin{prf}
It follows from Lemma \ref{thm: equiv-1st-eigenv} that $\lambda=\lambda(t,r)$ satisfies the equation
\begin{equation}\label{eqn:F-funct}
F(t,r;\lambda):=r\varphi(t,\lambda)\Psi(t,\lambda)
-\frac{1}{m(\lambda)}=0,\ (t,r)\in (0,+\infty)\times (0,\epsilon_0].
\end{equation}

Define
\begin{equation}\label{324}
u(\lambda):=\frac{\partial }{\partial \lambda}\left(\frac{1}{m(\lambda)}\right).
\end{equation}
Since $\lambda_1$ is simple pole of $m(\lambda)$ by Lemma \ref{weyl}$(iii)$, we can see that
$u(\lambda_1)\neq 0$. And from \eqref{small}, one can see that $u(\lambda(t,r))\neq 0$
for any $(t,r)\in (0,+\infty)\times (0,\epsilon_0]$, where $\epsilon_0$ is selected as same as the one in Theorem~\ref{thm:property-1st-eigenv-funct}.
By Remark~\ref{rem:bdd-r-lambda}, we can set
$$
\lambda_0:=\min\{\lambda(t,r):\ (t,r)\in[0,+\infty)\times[0,\epsilon_0]\}.
$$
For any $N\in \BN$, since $\varphi(x,\lambda)$ and $\Psi(x,\lambda)$ are both continuous differentiable on
$(x,\lambda)\in [0,N]\times[\lambda_0,\lambda_1]$,
we can choose a sufficiently small $\epsilon_N\in (0,\epsilon_0],$ such that
\begin{equation}\label{dev}
r\left|\frac{\partial[\varphi(x,\lambda)\Psi(x,\lambda)]}{\partial \lambda}\right|<\min\{|u(\lambda)|: \ \lambda\in[\lambda_0,\lambda_1]\}
\end{equation}
by Lemma \ref{weyl}$(iii)$.
As a result, for $(t,r)\in [0,N]\times[0,\epsilon_N]$ and $\lambda\in[\lambda_0,\lambda_1]$, one can see that
$$
\frac{\partial F(t,r,\lambda)}{\partial \lambda}=r\frac{\partial[\varphi(t,\lambda)\Psi(t,\lambda)]}{\partial \lambda}
-u(\lambda)\neq 0.
$$
According to the existence theorem for implicit functions, there exists unique implicit function $\lambda=\lambda(t,r)$ on $\in[0,N]\times [0,\epsilon_N]$
such that
$F(t,r,\lambda(t,r))=0$, and $\lambda(t,r)$ is continuous differentiable with the partial derivative formula that
\begin{equation}\label{eqn:diff-eqn-sys}
\left\{\aligned
&\frac{\partial \lambda(t,r)}{\partial t}\left\{u(\lambda)-r\frac{\partial[\varphi(t,\lambda)\Psi(t,\lambda)]}{\partial \lambda}\right\}
=r\frac{\partial[\varphi(t,\lambda)\Psi(t,\lambda)]}{\partial t},\\
&\frac{\partial \lambda(t,r)}{\partial r}\left\{u(\lambda)-r\frac{\partial[\varphi(t,\lambda)\Psi(t,\lambda)]}{\partial \lambda}\right\}
=\varphi(t,\lambda)\Psi(t,\lambda),
\endaligned\right.
\end{equation}
where $\lambda=\lambda(t,r)$. Since
$$\frac{\partial \varphi(x,\lambda)}{\partial x}, \
\frac{\partial \Psi(x,\lambda)}{\partial x},\ \frac{\partial[\varphi(x,\lambda)\Psi(x,\lambda)]}{\partial \lambda}\in AC_\BR[0,+\infty),
$$
it follows from \eqref{eqn:diff-eqn-sys} that
$$
\frac{\partial \lambda(t,r)}{\partial t}, \ \frac{\partial \lambda(t,r)}{\partial r}\in AC_\BR[0,N].
$$

According to \eqref{eqn:eigenf-perturb-probl}, the eigenfunction associated to $\lambda(t,r)$ satisfies
\begin{equation}\label{equ-sol}
\Phi(x,\lambda(t,r))=\left\{\aligned
&\varphi(x,\lambda(t,r)),\ &&x\in[0,t),\\
&c(t,r)\Psi(x,\lambda(t,r)),\  &&x\in(t,+\infty),
\endaligned\right.
\end{equation}
where $c(t,r)=\frac{\varphi(t,\lambda(t,r))}{\Psi(t,\lambda(t,r))}$.
Since $\varphi(t,\lambda(t,r))$ and $\Psi(t,\lambda(t,r))$ are both continuous differentiable on $(t,r)$ and
$\Psi(x,\lambda(t,r))>0$ for $x\in(0,+\infty)$, one can see that $c(t,r)$ is
continuous differentiable on $(t,r)\in(0,N)\times [0,\epsilon_N]$, and hence $\Phi(x,\lambda(t,r))$ is continuous
differentiable on $(t,r)\in(0,N)\times [0,\epsilon_N]$ and $t\neq x$.
We should note that for fixed $x\in(0,+\infty)$,
$$
\frac{\partial \Phi(x,\lambda(t,r))}{\partial r}=\left\{\aligned
&\frac{\partial\varphi}{\partial\lambda}\frac{\partial \lambda(t,r)}{\partial r},\ &&t\in(x,+\infty),\\
&\left[\frac{\partial c}{\partial\lambda}\Psi+c\frac{\partial\Psi}{\partial\lambda}\right]\frac{\partial \lambda(t,r)}{\partial r},\  &&t\in(0,x),
\endaligned\right.
$$
where $\varphi=\varphi(x,\lambda(t,r)), \Psi=\Psi(x,\lambda(t,r))$ and $c=c(t,r)$. One can verify that
$$
\frac{\partial \Phi(x,\lambda(x-0,r))}{\partial r}=\frac{\partial \Phi(x,\lambda(x+0,r))}{\partial r}=
\frac{\partial\varphi(x,\lambda(x,r))}{\partial\lambda}\cdot\frac{\partial \lambda(x,r)}{\partial r},
$$
which means $\Phi(x,\lambda(t,r))$ is continuous differentiable with respect to  $r\in[0,\epsilon_N]$.
This proves the  conclusion $(ii)$.

Finally, we prove the formula \eqref{eq:part-deriv-r}.

On one hand, under the hypothesis ${\bf (H_1)}$,
since $\Phi(x,\lambda(t,r))$ satisfies
\begin{equation}\label{sol}
-\Phi''+[q-r\delta(x-t)]\Phi=\lambda(t,r)\Phi, \ \Phi(0)=0,
\end{equation}
set $v=\frac{\partial \Phi}{\partial r}$, then $v$ satisfies
\begin{equation}\label{dev-sol}
-v''+[q-r\delta(x-t)]v-\delta(x-t)\Phi=\frac{\partial \lambda(r,t)}{\partial r} \Phi+\lambda(t,r) v,\ v(0)=0.
\end{equation}
Here we claim that
\begin{equation}\label{lim}
\lim\limits_{x\rightarrow +\infty} \Phi(x,\lambda(t,r))v'(x,\lambda(t,r))-\Phi'(x,\lambda(t,r))v(x,\lambda(t,r))=0.
\end{equation}
Indeed, set $u=\frac{\partial \chi}{\partial \lambda}$, where $\chi$ is defined as in Lemma \ref{weyl}$(i)$, then
 $\chi$ satisfies
\begin{equation}\label{equ-weyl}
-\chi''+q\chi= \lambda \chi,\ \chi(0,\lambda)=1,\ \chi'(0,\lambda)=m(\lambda),
\end{equation}
and $u$ satisfies
\begin{equation}\label{equ-u}
-u''+qu= \lambda u+\chi,\ u(0,\lambda)=0,\ u'(0,\lambda)=\int^{+\infty}_0 \chi^2(x,\lambda) \mathrm{d}x
\end{equation}
by Lemma \ref{weyl}$(ii)$.
Multiplying two sides of the equality \eqref{equ-u} and \eqref{equ-weyl} by $\chi$ and $u$, respectively and
integrating their difference over $[0,+\infty)$, one can verify that
\begin{equation}\label{equ-lim}
\lim\limits_{x\rightarrow +\infty} (-u'\chi+\chi' u)(x)=(-u'\chi+\chi' u)(0)+
\int^{+\infty}_0 \chi^2(x,\lambda) \mathrm{d}x=0.
\end{equation}
According to the expression \eqref{equ-sol} of $\Phi(x,\lambda(t,r))$ and \eqref{Psi} of $\Psi(x,\lambda(t,r))$, for sufficiently large $x$, one can verify
\begin{equation}\label{42}
\Phi v'-\Phi'v=
c^2(t,r)\left[\frac{\partial \Psi}{\partial \lambda} \Psi'- (\frac{\partial \Psi}{\partial \lambda})' \Psi \right]\frac{\partial \lambda}{\partial r}
=\frac{1}{m^2(\lambda(t,r))}\left[u \chi'- u' \chi\right]\frac{\partial \lambda}{\partial r},
\end{equation}
 then
the equation \eqref{lim} holds due to both of \eqref{42} and \eqref{equ-lim}.

Multiplying two sides of the equality \eqref{sol} and \eqref{dev-sol} by $v$ and $\Phi$, respectively and
integrating their difference over $[0,+\infty)$, one can verify that
\begin{equation}\label{eq:part-deriv-r-orgin}
\frac{\partial \lambda(t, r)}{\partial r}\int^{+\infty}_0 \Phi^2(x,\lambda(t,r)) d x=-\Phi^2(t,\lambda(t,r))
\end{equation}
on the basis of \eqref{lim}.
If  $\Phi$ is a normalized eigenfunction, then the above equality is just \eqref{eq:part-deriv-r}.

On the other hand, under the hypothesis  ${\bf (H_2)}$,
 since $\Phi(x,\lambda(t,r))$ satisfies
\begin{equation}
-\Phi''+[q-r\delta(x-t)]\Phi=\lambda(t,r)\Phi, \ \Phi(0)=0=[\Phi, f](+\infty)\sin\beta-[\Phi, g](+\infty)\cos\beta,
\end{equation}
set $v=\frac{\partial \Phi}{\partial r}$, then $v$ satisfies
\begin{equation}
-v''+[q-r\delta(x-t)]v-\delta(x-t)\Phi=\frac{\partial \lambda(r,t)}{\partial r} \Phi+\lambda(t,r) v,\ v(0)=0=[v, f](+\infty)\sin\beta-[v, g](+\infty)\cos\beta.
\end{equation}
Here we can also claim that \eqref{lim} holds, that is
\begin{equation}\label{equ-lim-lim-cir}
\lim\limits_{x\rightarrow +\infty} [\Phi,v](x,\lambda(t,r))=0,
\end{equation}
due to the real-valued nature of $\Phi, v$.
Indeed,
 $\chi$ satisfies
\begin{equation}
-\chi''+q\chi= \lambda \chi,\ \chi(0,\lambda)=1,\ [\chi, f](+\infty)\sin\beta-[\chi, g](+\infty)\cos\beta=0,
\end{equation}
where $\chi$ is defined in Lemma \ref{weyl}$(i)$.
Set $u=\frac{\partial \chi}{\partial \lambda}$. Then, since $f,g$ are independent of $\lambda$, on has that $u$ satisfies
\begin{equation}
-u''+qu= \lambda u+\chi,\ u(0,\lambda)=0=[u, f](+\infty)\sin\beta-[u, g](+\infty)\cos\beta.
\end{equation}
Since both of $\chi$ and $u$ belong to $\mathcal{D}_{max}$, we have
\begin{equation}
\lim\limits_{x\rightarrow +\infty} [\chi,u](x,\lambda(t,r))= \lim\limits_{x\rightarrow +\infty} [\chi,f][u,g]-[u,g][\chi,f]=0
\end{equation}
by the bracket decomposition property (see \cite[Lemma 10.4.2]{Zettl2005}).
Then the formula \eqref{equ-lim-lim-cir} holds according to \eqref{42}.
Similar to the case under ${\bf (H_1)}$, we can also get the equation \eqref{eq:part-deriv-r-orgin}, which gives the formula \eqref{eq:part-deriv-r} under the hypothesis ${\bf (H_2)}$. The proof is completed.                       
\end{prf}

\bigskip
\section{Recovery of potential by the first eigenvalue function}\label{sec:potential-recovery}
\medskip

In this section, we show that the unique concrete expression of the potential $q$ is given by the first eigenvalue function of $({\bf E}_q)$.

\begin{thm}\label{thm:main-thm}
Let $\lambda(t,r)$ be the first eigenvalue function of $({\bf E}_q)$. Under either of the hypotheses ${\bf (H_1)}$ and ${\bf (H_2)}$, the potential $q(x)$ is determined by $\lambda(t,r)$ uniquely, and can be expressed as
\begin{equation}\label{eqn:reconstruct-formula}
q(x)=\frac{\varphi''_0(x)}{\varphi_0(x)}+\lambda_1,\
\varphi_0(x):=\sqrt{-\frac{\partial \lambda(x,0)}{\partial r}},\ x\in(0,+\infty),
\end{equation}
where $\lambda_1$, which equals to $\lambda(t,0)$, is the first eigenvalue of $({\bf E}_q)$.
\end{thm}

\begin{prf} By Lemma~\ref{lem:dis-spc-Str-lp} and Lemma~\ref{lem:dis-spc-Str-lc}, we can know that under either of the hypotheses ${\bf (H_1)}$ and ${\bf (H_2)}$, the first eigenvalue function $\lambda(t,r)$ of $({\bf E}_q)$ exists.
Letting $r\rightarrow 0^{+}$ in \eqref{eq:part-deriv-r} and using Theorem~\ref{thm:diff-1st-eigenv-funct}$(ii)$,  we have
\begin{equation}\label{41}
\frac{\partial \lambda(t,0)}{\partial r}=-\Phi_0^2(t,\lambda_1),
\end{equation}
and hence the desired $\varphi_0(x):=\sqrt{-\frac{\partial \lambda(x,0)}{\partial r}}$ is a non sign-changed normalized solution of $$-\varphi_0''(x)+q(x)\varphi_0(x)=\lambda_1\varphi_0(x),\ x\in(0,+\infty),$$
which gives the formula \eqref{eqn:reconstruct-formula}.
\end{prf}

\begin{rem}
According to Lemma \ref{weyl}$(ii)$ and the definition of $\Psi$ in \eqref{Psi}, one can verify
$$\frac{\partial }{\partial \lambda}\left(\frac{1}{m(\lambda(t,r))}\right)=
-\int^{+\infty}_0 \Psi^2(x,\lambda) \mathrm{d}x.$$
Letting $r\to0^+$ in the second equality in \eqref{eqn:diff-eqn-sys}, one can also directly get \eqref{41}.
\end{rem}

\begin{rem}\label{rem:approx-recov}
$(i)$ The above theorem also implies the uniqueness:  if $\lambda(t,r;q_1)=\lambda(t,r;q_2)$, then $q_1=q_2$.

\smallskip
\noindent
$(ii)$ Note that $\frac{\partial\lambda(t,0)}{\partial r}$ is a right partial derivative at $r=0$.
In fact, if $\lambda(t,r)$ is the first eigenvalue function of $({\bf E}_q)$ with a unknown potential $q$, then, to recover $q$, we only needs to know the values of
$\lambda(t,r_n)$ for $r_n\to 0^+$ as $n\to +\infty$ and $t\in (0,+\infty)$, since one has
$$\frac{\partial \lambda(t,0)}{\partial r}=\lim\limits_{n\to +\infty}\frac{\lambda(t,r_n)-\lambda(t,0)}{r_n}.$$
\end{rem}

\begin{rem}\label{rem: other-case}
If $({\bf \tilde{E}}_q)$ does not have eigenvalue, or has oscillatory eigenvalue (i.e., every eigenfunction has an infinite number of zeros), then we can not recover the unknown potential through our method in the present paper.
For instance,

\mednoind
$(i)$ if $\tau_q y$ in $({\bf \tilde{E}}_q)$ is in the limit-circle case and oscillatory at $+\infty$, then the spectrum is unbounded above and below, and every eigenfunction has infinite number of zeros in $(0,+\infty)$, and hence there is no principle eigenvalue (see \cite[Theorem 10.12.1(7)]{Zettl2005});

\mednoind
$(ii)$ when $\tau_q y$ in $({\bf \tilde{E}}_q)$ is in the limit-point case at $+\infty$, even if the spectrum is discrete, then the absence of boundedness below of the spectrum still will cause the non-existence of non-oscillatory eigenvalue (see \cite[Theorem 10.12.1(8)(ii)]{Zettl2005}), namely, in the limit-point case such as in our paper, the boundedness below of the spectrum is necessary;

\mednoind
$(iii)$ if $\tau_q y$ in $({\bf \tilde{E}}_q)$ is in the limit-point case at $+\infty$, and $q$ is bounded below with bound $q_0$, there exists essential spectrum. For the eigenvalue $\lambda>q_0$, the eigenfunction corresponding to $\lambda$ has an infinite numbers of zeros (see \cite[Theorem 10.12.1(8)(iii)]{Zettl2005}).  Therefore, in this case, we are only allowed to consider the eigenvalue below the essential spectrum, such as Lemma~\ref{lem:principle-eigenv-exist-case}$(i)$.

\end{rem}

To sum up, Theorem \ref{thm:main-thm} only gives the unique recover of potential by the first eigenvalue function. However, there still remains many unsolved questions, such as the existence of the first eigenvalue function, stability of reconstructing the potential, and so on. So, we can also propose many similar motivating questions as the ones at the end of \cite{ZQC2024}.

\bigskip

\section*{Acknowledgement}
\smallskip

This research was partially supported by the National Natural Science Foundation of China [Grant numbers 12271299, 12071254, and 11701327] and Shandong Provincial Fund [ZR2024MA005]. The authors are grateful to Professor Xiaoping Yuan for his helpful discussions and guidance.

\end{document}